\documentclass[12pt,fleqn]{article}
\usepackage{graphicx}

\usepackage{latexsym}

\usepackage{amsmath}
\usepackage{amsthm}
\usepackage{amssymb}
\usepackage{amsfonts}

\begin{document}

\newcommand{\rf}[1]{(\ref{#1})}
\newcommand{\rff}[2]{(\ref{#1}\ref{#2})}

\newcommand{\ba}{\begin{array}}
\newcommand{\ea}{\end{array}}

\newcommand{\be}{\begin{equation}}
\newcommand{\ee}{\end{equation}}

\newcommand{\const}{{\rm const}}
\newcommand{\ep}{\varepsilon}
\newcommand{\Cl}{{\cal C}}
\newcommand{\rr}{\vec r}
\newcommand{\ph}{\varphi}

\newcommand{\e}{{\bf e}}

\newcommand{\m}{\left( \ba{r}}
\newcommand{\ema}{\ea \right)}
\newcommand{\mm}{\left( \ba{cc}}
\newcommand{\miv}{\left( \ba{cccc}}

\newcommand{\scal}[2]{\mbox{$\langle #1 \! \mid #2 \rangle $}}
\newcommand{\ods}{\par \vspace{0.5cm} \par}
\newcommand{\dis}{\displaystyle }
\newcommand{\mc}{\multicolumn}

\newtheorem{prop}{Proposition}
\newtheorem{Th}{Theorem}
\newtheorem{lem}{Lemma}
\newtheorem{rem}{Remark}
\newtheorem{cor}{Corollary}
\newtheorem{Def}{Definition}
\newtheorem{open}{Open problem}
\newtheorem{ex}{Example}
\newtheorem{exer}{Exercise}

\title{\bf How to improve the accuracy of \\ the discrete gradient method \\ in the one-dimensional case}

\author{
 {\bf Jan L.\ Cie\'sli\'nski}\thanks{\footnotesize
 e-mail: \tt janek\,@\,alpha.uwb.edu.pl}
\\ {\footnotesize Uniwersytet w Bia{\l}ymstoku,
Wydzia{\l} Fizyki, ul.\ Lipowa 41, 15-424
Bia{\l}ystok, Poland}
\\ {\bf Bogus{\l}aw Ratkiewicz}\thanks{\footnotesize
e-mail: \tt bograt\,@\,poczta.onet.pl} \thanks{\footnotesize permanent address: 
 I L.O.,
ul.\ \'Sr\'odmie\'scie 31, 16-300 August\'ow, Poland;}
\\ {\footnotesize Doctoral Studies, Wydzia\l \ Fizyki, Uniwersytet Adama Mickiewicza, Pozna\'n, Poland}
}

\date{}

\maketitle

\begin{abstract}
We present a new numerical scheme for one dimensional dynamical systems. This is a modification of the discrete gradient method and keeps its advantages, including the stability and the conservation of the energy integral. However, its accuracy is higher by several orders of magnitude. 
\end{abstract}

\ods




\pagebreak

The discrete gradient methods were introduced many years ago in order to integrate numerically $N$-body systems of classical mechanics with possible applications in molecular dynamics and celestial mechanics \cite{LaG} (see also \cite{HLW}). In this paper we consider the one-dimensional case:
\be  \label{Newton}
\dot p =  - V'(x) \ , \quad \dot x = p \ , 
\ee
where $V(x)$ is a potential, and dot and prime denote differentiation with respect to $t$ and $x$, respectively. In this case discrete gradient methods reduce to the so called modified midpoint rule: 
\be \ba{l} \label{dis-grad} \dis
\frac{x_{n+1} - x_n}{\ep} = \frac{1}{2} \left( p_{n+1} + p_n \right) \ . \\[3ex] \dis
\frac{p_{n+1} - p_n}{\ep} =  - \frac{ V (x_{n+1}) - 
V (x_n) }{x_{n+1} - x_n} \ . 
\ea \ee
One can easily prove that the total energy is preserved by this scheme:
\be  \label{energy}
 \frac{1}{2} p_n^2 + V (x_n) = E = \const \ , 
\ee
The modified midpoint rule has been extended, in a natural way, on the three-dimensional case and on systems of particles, exactly preserving the total energy, the total linear momentum, and the total angular momentum of the system \cite{LaG}. 

More recently discrete gradient methods have been developed in the context of geometric numerical integration \cite{MQ}, see \cite{IA,Gon}.   In particular, Quispel and his collaborators constructed numerical integrators preserving integrals of motion of a given system of ordinary differential equations \cite{QC,QT,MQR2}.  

In general, geometric numerical integrators are very good in preserving qualitative features of simulated differential equations but it is not easy to enhance their accuracy. Our research is concentrated on improving  the efficiency of the discrete gradient method without loosing its outstanding qualitative advantages. 

\ods

In a previous paper we compared several discretizations of the simple pendulum equation with a special stress on the long-time behaviour \cite{CR-long}. The discrete gradient scheme was among the best ones, especially when large energies (rotational motion) and the neighbourhood of the separatrix were concerned. 
In the paper \cite{CR-long} we proposed a modification of the discrete gradient scheme \rf{dis-grad}. Assuming the stable equilibrium at $x=0$, we replaced $\ep$ by the function $\delta_0 = \delta_0 (\ep)$
\be
\delta_0 = \frac{2}{\omega_0} \tan \frac{\omega_0\ep}{2} \ ,
\ee
where $\omega_0 = \sqrt{V''(0)}$. The motivation was to preserve (almost exactly) small oscillations around $x =0$, when the pendulum can be treated as a harmonic oscillator.  
The classical harmonic oscillator admits the so called {\it exact} (or the best) discretization \cite{Ag,CR-ade}. The exact discretization of the harmonic oscillator equation was earlier used to derive an orbit-preserving discretization of the classical Kepler problem \cite{Ci-Kep}. 
The result obtained in the paper \cite{CR-long} has been quite satisfying: in the case of small oscillations our method was better by 4 orders of magnitude than all other considered schemes (including the discrete gradient method). In the case of other initial conditions the new method was comparable with the discrete gradient method (by the way: both methods are of second order). 

\ods

\begin{figure}
\caption{\small $\delta_n/\ep$ \ for $\ep = 0.02$ (dashed line) and $\ep = 0.1$ (solid line).}
 \label{delta}  \par
\includegraphics[height=0.46\textheight]{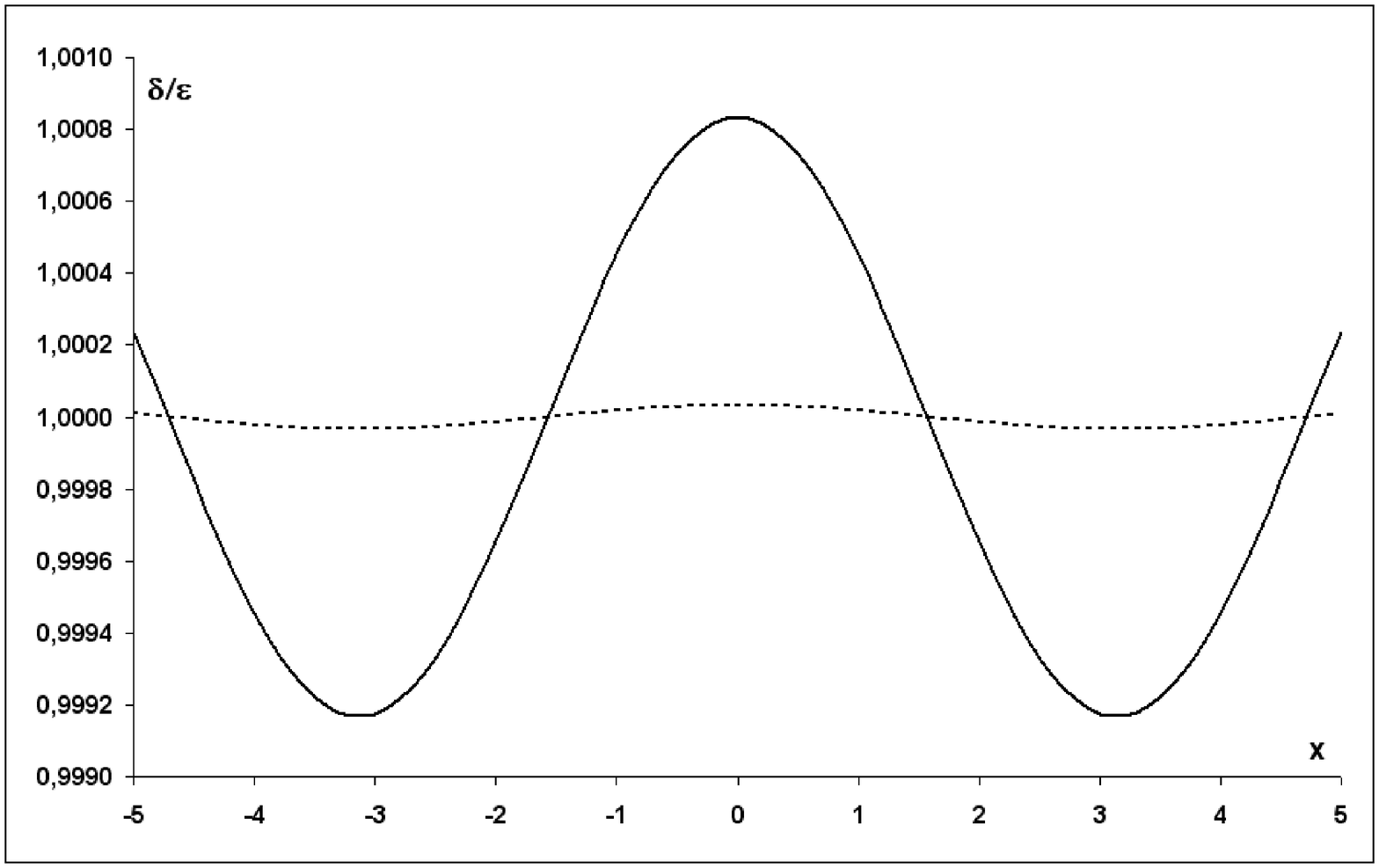} \par
\end{figure}

As the main result of this paper we propose another, much more efficient,  
modification of the discrete gradient scheme:
\be \ba{l} \label{corr} \dis
\frac{x_{n+1} - x_n}{\delta_n} = \frac{1}{2} \left( p_{n+1} + p_n \right) \ . \\[3ex] \dis
\frac{p_{n+1} - p_n}{\delta_n} =  - \frac{ V (x_{n+1}) - 
V (x_n) }{x_{n+1} - x_n} \ ,
\ea \ee
where $\delta_n$ is a function defined by
\be \ba{l}  \label{deltan} \dis
\delta_n =
\frac{2}{ \omega_n } \tan\frac{ \ep \omega_n  }{2} \ , 
\quad ({\rm if } \  V'' (x_n) > 0) \ ,  \\[3ex]
\delta_n = \ep \ ,  
\quad ({\rm if } \  V'' (x_n) = 0) \ , \\[2ex] \dis
\delta_n =
\frac{2}{ \omega_n } \tanh \frac{ \ep \omega_n }{2} \ , 
\quad ({\rm if } \  V'' (x_n) < 0) \ ,
\ea \ee
$\ep$ denotes the time step, and, finally, 
\be \label{omegan}
\omega_n = \sqrt{  | V'' (x_n) | } \ .  
\ee
In fact we have replaced $\ep$, appearing in formulas \rf{dis-grad}, by a variable $\delta_n$ depending not only on $\ep$ but also on $x_n$. The crucial point is to choose the best form of the function $\delta_n$. In order to obtain \rf{deltan} we have to require that 
the modified scheme \rf{corr} is {\it locally exact} \cite{Ci-focm}. 
It means that the linearization of the scheme \rf{corr} 
around any fixed $x$, i.e., 
\be \ba{l} \label{corr-lin} \dis
\frac{\xi_{n+1} - \xi_n}{\delta_n} = \frac{1}{2} \left( p_{n+1} + p_n \right) \ . \\[3ex] \dis
\frac{p_{n+1} - p_n}{\delta_n} =  - V'(x) - \frac{1}{2} (\xi_{n+1} + \xi_n) V''(x) \ ,
\ea \ee
(where $\xi_n = x_n - x$) coincides with the {\it exact discretization} of the linearization of the considered system \rf{Newton}  
\be  \label{Newtonlin}
 \frac{d p}{d t} = - V'(x) - V''(x) \xi \ , \qquad 
\frac{d \xi}{d t} = p \ , 
\ee
where $x$ is treated as a constant 
(equations \rf{Newtonlin}, equivalent to the harmonic oscillator with a constant force, admit the exact discretization).  
Thus we can express $\delta_n$ in terms of $V''(x)$. Then, identifying $x$ with $x_n$, we arrive at formulas \rf{deltan}. The detailed derivation of the numerical scheme \rf{corr}  will be presented elsewhere \cite{Ci-focm}. 
We point out that $\delta_n$ is a function which is almost constant 
(especially for small $\ep$), namely: $\delta_n \approx \ep$, see Fig.~\ref{delta}.    
\ods
The replacement $\ep \rightarrow \delta_n$  works very well for the first order system \rf{corr}. However, if we try to replace $\ep$ by $\delta_n$   in the second order discrete equation for $x_n$ (see Eq.~(46) in the paper \cite{CR-long}), then we get a difference scheme which is only a little bit better than the scheme with constant $\delta = \delta_0$.  

 \ods
In practical implementation we use the implicit scheme \rf{corr} as 
the corrector, while taking as the predictor the explicit scheme 
\be  \ba{l} \label{pred}  \dis
  x_{n+1} = x_n + 
\frac{\sin (\omega_n \ep) }{\omega_n } p_n
- \frac{1 - \cos(\omega_n\ep)}{\omega_n^2} \
V'(x_n) \ , \quad ({\rm if } \  V'' (x_n) > 0) \ ,\\[2ex] \dis
x_{n+1} = x_n  +
\ep p_n - \frac{1}{2} \ep^2 \ V' (x_n) \ , \quad ({\rm if } \  V'' (x_n) = 0) \ ,\\[2ex] \dis
x_{n+1} = x_n  +
\frac{\sinh (\omega_n \ep) }{\omega_n } p_n
- \frac{1 - \cosh (\omega_n\ep)}{\omega_n^2} \
V' (x_n) \ , \ ({\rm if } \  V'' (x_n) < 0)  \ .
\ea \ee
In order to obtain \rf{pred} one has to eliminate $p_{n+1}$ from the system \rf{corr} and expand the result in the Taylor series with respect to $x_{n+1}- x_n$, \ leaving only linear terms. 

The numerical scheme \rf{pred} has the same order (third) as \rf{corr}. However, this is of some advantage only for very small $\ep$ and for very short times (thus \rf{pred} can serve as a very good predictor). As far as the long-time behaviour is concerned the scheme \rf{pred} is not good and yields solutions with wrong qualitative behaviour.   

\ods
The proposed numerical integrator \rf{corr} (we will refer to it as {\it locally exact  discrete gradient method}) has important advantages: 
\begin{itemize}
\item exact conservation of the energy integral (i.e., eq.\ \rf{energy} holds),
\item higher order (third) as compared with the discrete gradient method, 
\item high stability and accuracy,
\item very good long-time behaviour of numerical solutions.
\end{itemize} 

\ods
The accuracy of our new numerical scheme was tested on the case of the simple pendulum equation ($V (x) = - \cos x$). We compared it with the standard leap-frog scheme, the discrete gradient method and modified discrete gradient method (introduced in \cite{CR-long}). 
All motions of the simple pendulum are periodic, so we focus our attention on the relative error of the period of considered discretizations. The details of our numerical experiments are explained in \cite{CR-long}.  Here we just mention that for simplicity we assume $x_0 = 0$. 
In this case $p_0^2 = 2 E + 2$.  

\begin{figure}
\caption{\small Relative error of the period as a function of $p_0$ for $\ep = 0.02$. 
White triangles: leap-frog, white diamonds: discrete gradient, black diamonds: modified discrete gradient ($\delta = \const$), black squares: locally exact discrete gradient. }
 \label{e002}  \par
\includegraphics[height=0.40\textheight]{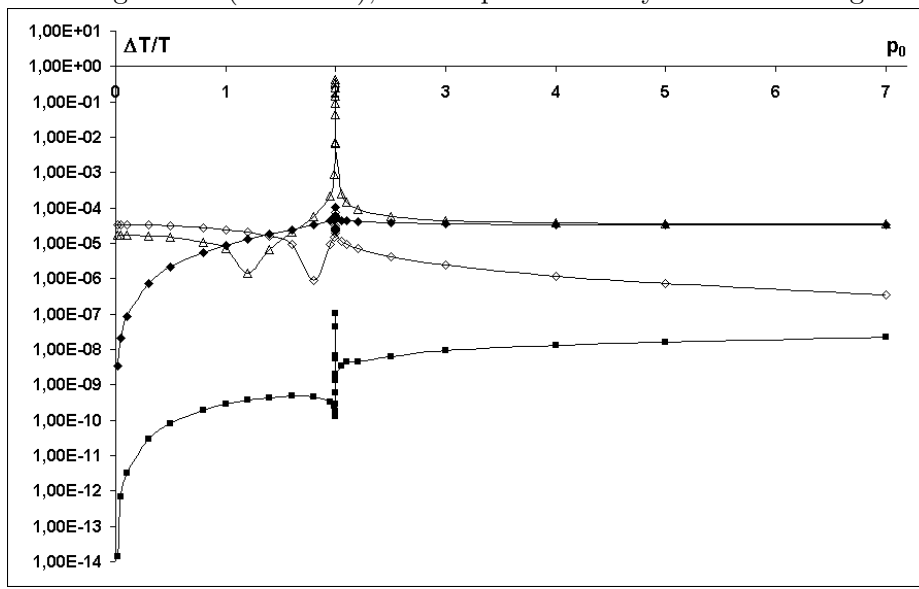} \par
\end{figure}

Both discrete gradient methods studied in \cite{CR-long} are very stable and our new method shares this property as well. 
The accuracy of the new method is surprisingly high, especially for small 
(but not necessarily very small) time-steps. 
As an example we consider the case $\ep=0.02$ (see Fig.~\ref{e002}). 
For small oscillations (e.g., $p_0 = 0.02$) the accuracy of our new method  is  greater by 5 orders of magnitude as compared to the modified discrete gradient method, and 9 (nine!) orders of magnitude better than leap-frog or discrete gradient method. Actually, the locally exact discrete gradient method is much more accurate than any other considered method for any initial conditions, see Fig.~\ref{e002} (note that the scale on the vertical axis is logarithmic).  Usually other methods are worse at least by 4 orders of magnitude, except the discrete gradient method which, for large $p_0$, is worse ``only'' by 2 orders of magnitude and this difference seems to diminish for larger energies.  

\begin{figure}
\caption{\small Relative error of the period as a function of $\ep$ for 
$p_0 = 0.02$. Symbols: the same as in figure~\ref{e002}. } 
 \label{p002}  \par
\includegraphics[height=0.40\textheight]{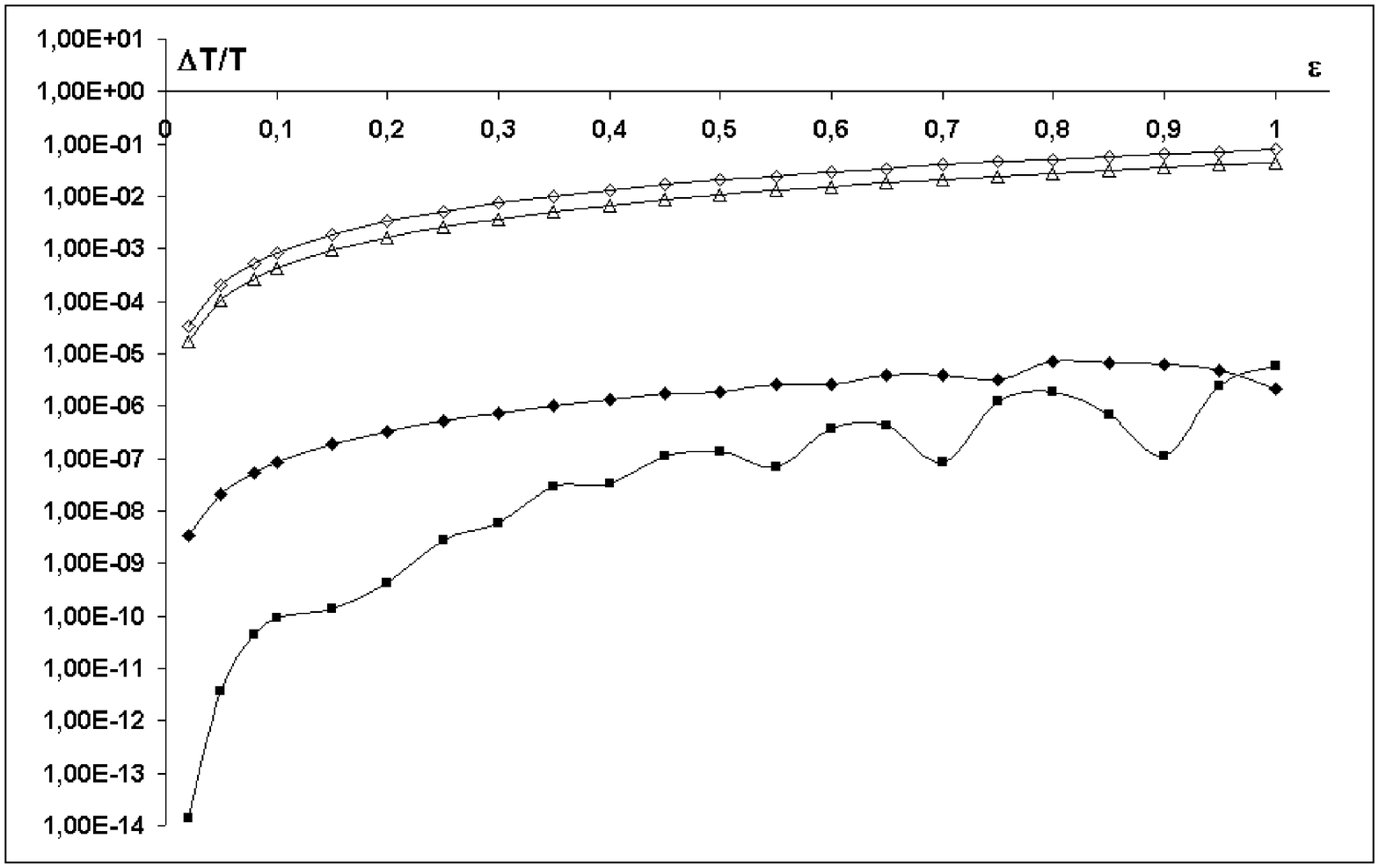} \par
\end{figure}

\begin{figure}
\caption{\small Relative error of the period as a function of $\ep$ for 
$p_0 = 1.21$ (``resonance value'' 
for the leap-frog scheme). Symbols: the same as in figure~\ref{e002}.  } 
 \label{p121}  \par
\includegraphics[height=0.40\textheight]{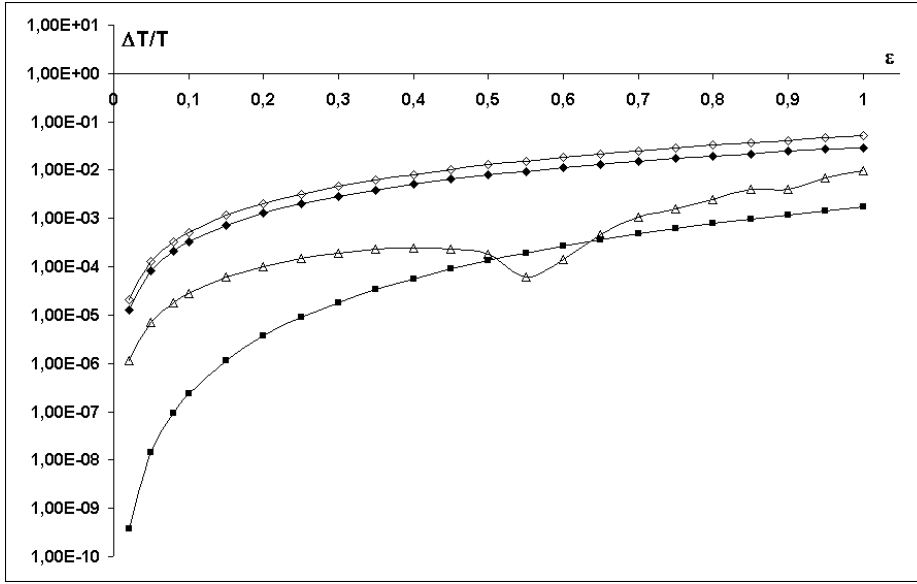} \par
\end{figure}

\begin{figure}
\caption{\small Relative error of the period as a function of $\ep$ for $p_0 = 1.8$. Symbols: the same as in figure~\ref{e002}. } 
 \label{p18}  \par
\includegraphics[height=0.40\textheight]{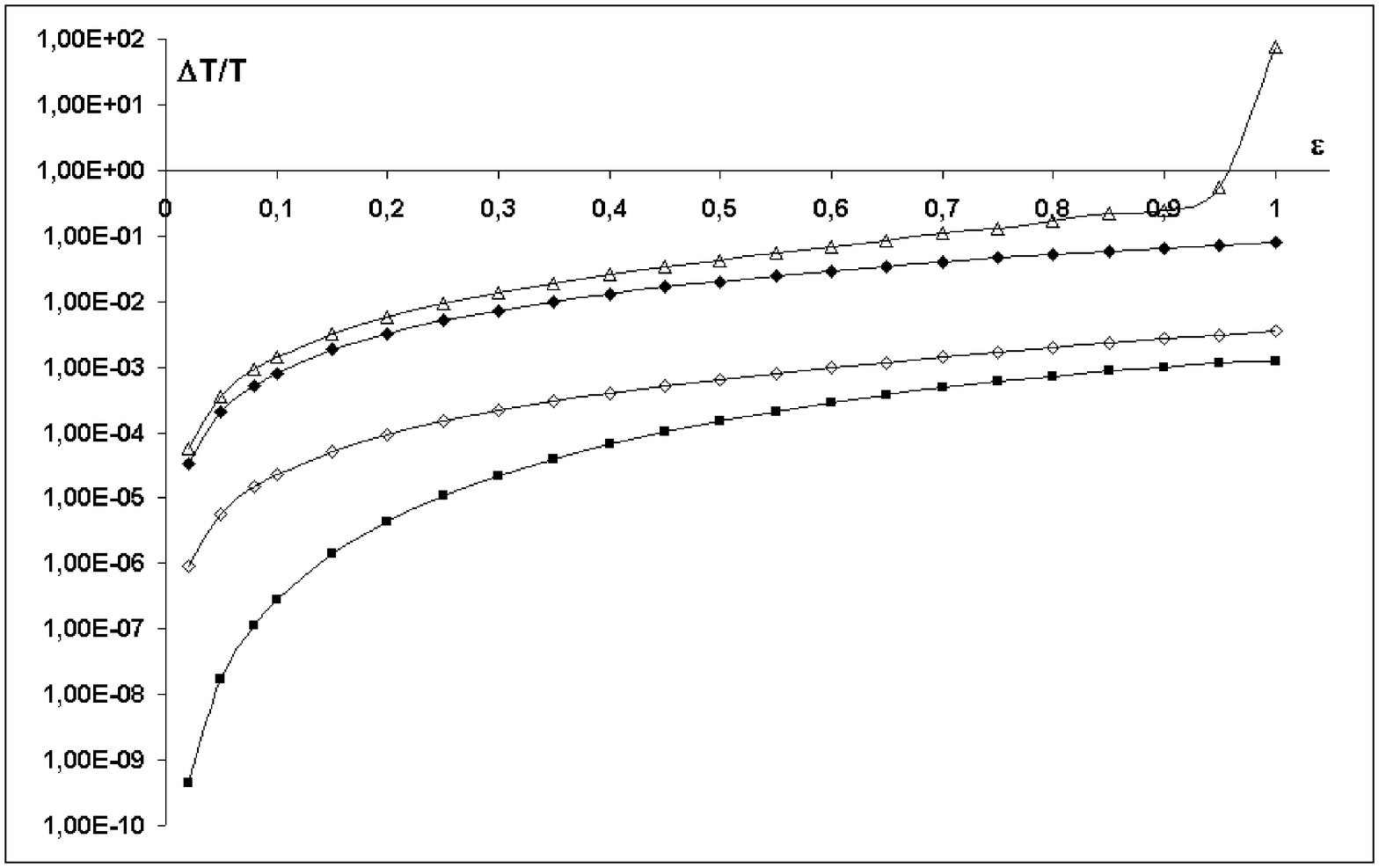} \par
\end{figure}

We present also the time-step dependence of the considered numerical schemes for three chosen energies (see Fig.~\ref{p002}, Fig.~\ref{p121} and Fig.~\ref{p18}, note that these graphs are semilogarithmic). The locally exact discrete gradient scheme is almost always the best with very few  exceptions. Namely, in the case of small oscillations and larger time-steps the modified discrete gradient scheme yields similar results, see Fig.~\ref{p002}. Then, for larger amplitudes and larger time-steps the discrete gradient scheme is not much worse, see Fig.~\ref{p18}. Finally,  the leap-frog 
is comparable with our new method when $p_0$ is close to a peculiar value $p_0 \approx 1.21$, see Fig.~\ref{p121}. This is the somewhat mysterious ``resonance value'' of $p_0$ for the leap-frog scheme, already reported in \cite {CR-long}.  For this initial condition the leap-frog method is exceptionally accurate and, for larger time-steps, is comparable even 
with the locally exact discrete gradient scheme.

\begin{figure}
\caption{\small $x_n$ as a function of $n$, very near the separatrix ($p_0 = 1.9999999999$), for $\ep = 0.9$. Symbols: the same as in figure~\ref{e002}. The solid line corresponds to the exact (continuous) solution.    } 
 \label{pod-sep}  \par
\includegraphics[height=0.40\textheight]{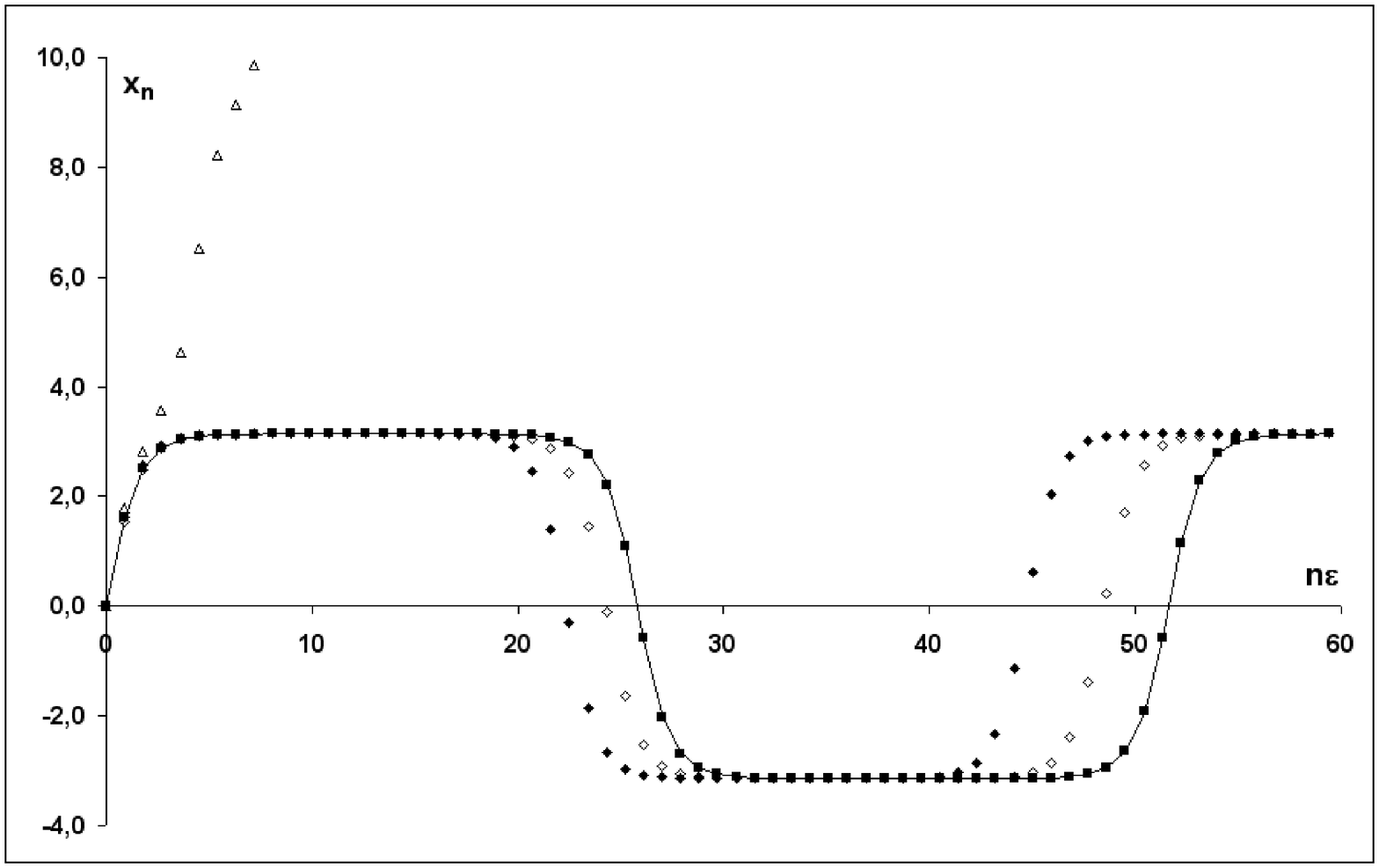} \par
\end{figure}

\begin{figure}
\caption{\small $x_n$ as a function of $n$, very near the separatrix ($p_0 = 2.0000000001$), for $\ep = 0.7$. Symbols: the same as in figure~\ref{e002}. The solid line corresponds to the exact (continuous) solution.    } 
 \label{nad-sep}  \par
\includegraphics[height=0.40\textheight]{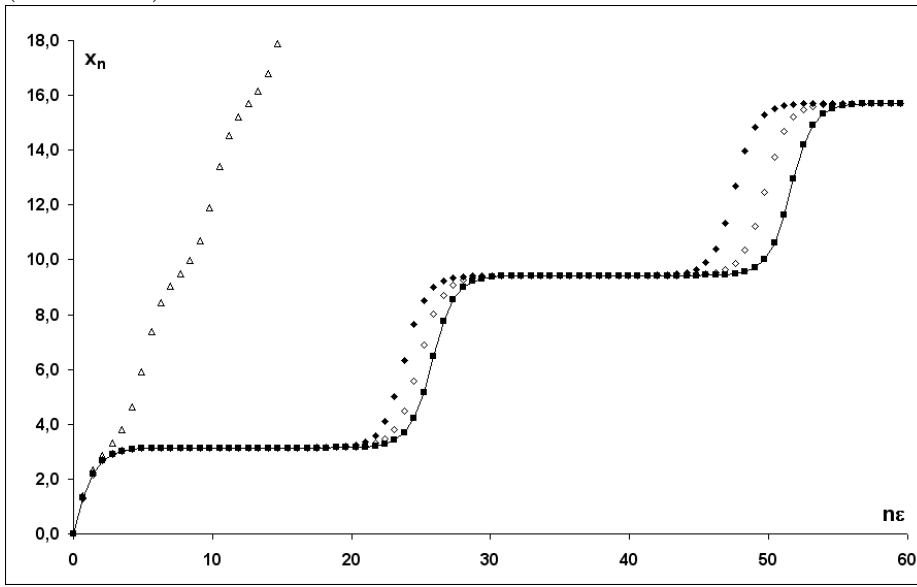} \par
\end{figure}

The neighbourhood of the separatrix ($p_0 \approx 2$) is most difficult to be simulated numerically. The discrete gradient method was relatively good at this region, see \cite{CR-long}. The locally exact discrete gradient method is excellent even in that case, see Fig.~\ref{pod-sep} and Fig.~\ref{nad-sep}. We are very close to the separatrix ($|p_0 - 2| = 10^{-10}$), $\ep$ is large, and nevertheless our new method simulates very accurately the motion of the pendulum. Discrete points $x_n$ practically lie on the continuous curve of the exact solution. The other two discrete gradient methods yield good results (at least qualitatively) while the leap-frog scheme fails to reproduce even the qualitative behaviour.  

\ods

The method presented in this paper can be extended and generalized 
on multidimensional cases and some other numerical integrators (including the implicit midpoint rule) \cite{Ci-focm}. 
The main idea of this novel approach consists in modifying numerical integrators by replacing $\ep$ by appropriate functions (e.g.,  $\delta_n$), depending on $\ep$ and independent variables, in order to obtain ``locally exact'' integrators (for instance, integrators which are locally equivalent to the exact discrete harmonic oscillator). Note that the time-step of the locally exact discrete gradient scheme \rf{corr} is equal to $\ep$ and is assumed to be constant. However, there are no obstacles to use the variable time-step and it is worthwhile to examinate this possibility in the future.

\pagebreak

\end{document}